\documentclass[a4paper, 12pt]{amsart}

\usepackage {algorithmic, algorithm, epsfig,nicefrac,amsthm}

\usepackage[utf8]{inputenc}

\usepackage{mathptmx}

\newtheorem{lemma}{Lemma}[section] \newtheorem{corollary}[lemma]{Corollary} \newtheorem{theorem}[lemma]{Theorem}

\theoremstyle{definition} \newtheorem{definition}[lemma]{Definition} \newtheorem{remark}[lemma]{Remark}

\def\ZZ{{\mathbb Z}}
\def\RR{{\mathbb R}}

\DeclareMathOperator{\rem}{rem}

\DeclareMathOperator{\Hilb}{Hilb}
\DeclareMathOperator{\ld}{ld}

\DeclareMathOperator{\parf}{par}

\def\ln{\operatorname{log}}

\begin {document}

\title[Unimodular triangulations of simplicial cones]{Unimodular triangulations of simplicial cones\\ by short vectors}

\author{Michael von Thaden}
\address{DekaBank,
	Deutsche Girozentrale,
	60325 Frankfurt am Main,
	Germany}
\email{\texttt{v\_thaden@t-online.de}}

\author{Winfried Bruns}
\address{Universit\"at Osnabr\"uck, Institut f\"ur Mathematik, 49069 Osnabr\"uck, Germany}
\email{wbruns@uos.de}



\maketitle

\begin{abstract}
We establish a bound for the length of vectors involved in a unimodular triangulation of simplicial cones. The bound is exponential in the square of the logarithm of the multiplicity, and improves previous bounds significantly. The proof is based on a successive reduction of the highest prime divisor of the multiplicity and uses the prime number theorem to control the length of the subdividing vectors.
\end{abstract}

\marginparwidth=5em
\marginparpush=10pt

\section{Introduction}

In this paper we discuss the triangulation of simplicial cones $C$ into unimodular subcones by ``short'' vectors. Length is measured by the basic simplex $\Delta_C$ of $C$ that is spanned by the origin and the extreme integral generators of $C$: we want to find an upper bound for the dilatation factor $c$ for which all subdividing vectors are contained in $c\Delta_C$. Roughly speaking, the larger the multiplicity $\mu$ of $C$ (given by the lattice normalized volume of $\Delta_C$), the more subdivision steps are to be expected, and they inevitably increase the length of the subdividing vectors. Therefore $\mu$ is the natural parameter on which estimates for $c$ must be based, at least for fixed dimension $d$.

A prominent case in which bounds for $c$ come up is the desingularization of toric varieties. The standard argument applied in this situation leads to rather bad bounds. A slight improvement was reached by Bruns and Gubeladze \cite[Theorem 4.1]{BG02a} who gave a bound that is better, but still exponential in $\mu$. The main result of this paper is a bound essentially of order $\mu^{\log \mu}$ for fixed dimension $d$ (Corollary \ref{simpler}). The next goal would be a bound that is polynomial in $\mu$, but we do not know if such exists.

It seems that the only general technique for triangulating a (simplicial) cone into unimodular subcones is successive stellar subdivision: one chooses an integral vector $x$ in $C$ and replaces $C$ by the collection of subcones that are spanned by $x$ and the facets of $C$ that are visible from $C$. This simple procedure allows successive refinement of triangulations if simultaneously applied to all cones that contain $x$.

We start from the basic and rather easy observation that a unimodular triangulation by iterated subdivision can be reached very quickly by short vectors if $\mu$ is a power of $2$. In order to exploit this observation for arbitrary $\mu$, two crucial new ideas are used: (i) not to diminish the multiplicity $\mu$ in every subdivision step (as usually), but to allow it to grow towards a power of $2$, and (ii) to control this process by the prime number theorem. While we formulate our results only for simplicial cones, they can easily be generalized (see Remark \ref{general}.) 

The bound in \cite[Theorem 4.1]{BG02a} was established in order to prove that multiples $cP$ of lattice polytopes $P$ can be covered by unimodular simplices as soon as $c$ exceeds a threshold that depends only on the dimension $d$ (and not on $P$ or $\mu$). If only unimodular covering is aimed at (for polytopes or cones), one can do much better than for triangulations: the threshold for $c$ has at most the order of $d^6$. In particular, it is independent of $\mu$. See Bruns and Gubeladze \cite[Theorem 3.23]{BG09}. A polynomial bound of similar magnitude for the unimodular covering of cones is given in \cite[Theorem 3.24]{BG09}. These polynomial bounds are based on the first part of von Thaden's PhD thesis \cite{vT07}, whereas the results of this paper cover the second part of \cite{vT07}. 

The most challenging problem in the area of this paper is to show that the multiples $cP$ of a lattice polytope $P$ have unimodular triangulations for \emph{all} $c\gg 0$. The best known result in arbitrary dimension is the Knudsen-Mumford-Waterman theorem that guarantees the existence of such $c$. We refer the reader to \cite[Chapter 3]{BG09} and to \cite{HPPS} for an up-to-date survey. The paper \cite{HPPS} contains an explicit upper bound for $c$.

For unexplained terminology and notation we refer the reader to \cite{BG09}.

\section{Auxiliary results}


Our first theorem will show that there is a sublinear bound in $\mu(C)$ on the length of the subdividing vectors in a unimodular triangulation if the multiplicity $\mu(C)$ of the cone $C$ is a power of $2$. We always assume that a simplicial cone $C$ is generated by its extreme integral generators, i.e., the primitive integral vectors contained in the extreme rays. The multiplicity then is the lattice normalized volume of the simplex $\Delta_C$ spanned by them and the origin.

Note that for unimodular cones $D$ the Hilbert basis $\Hilb(D)$ that appears in the theorem consists only of the extreme integral generators.

\begin{theorem}\label{power2}
Let $d \geq 3$ and let
$C = \RR_+ v_1+ \cdots + \RR_+ v_d \subset \RR^d$
be a simplicial $d$-cone with $\mu ( C ) = 2^l \ ( l \in \mathbb{N} )$.
Then there exists a unimodular triangulation $C = C_1 \cup \ldots \cup C_k $ such that
$$
\Hilb(C_j ) \subset \left( \frac{d}{2} \left( \frac{3}{2}\right)^l \right)\Delta_C,
\quad 1 \le j \le k.
$$
\end{theorem}

\begin{proof} The proof of this theorem is similar to the proof of Theorem 4.1 in \cite{BG02a}.
We consider the following sequence:
$$
h_k = 1, \quad k \leq 0,\qquad
h_k = \frac{1}{2}( h_{k-1} + \cdots + h_{k-d}), \quad k\geq 1.
$$
\noindent Because
$$
h_k - h_{k-1} = \frac{1}{2} h_{k -1} - \frac{1}{2} h_{k-d-1}
$$
for $k\geq 2$ and
$h_1 > h_l $ for $l \leq 0$, it follows by induction that this sequence is increasing.
Since for $k\geq 2$
$$
h_k = \frac{1}{2} h_{k-1} + \frac{1}{2} ( h_{k-2} + \cdots + h_{k-d-1}) - \frac{1}{2} h_{k-d-1} =
\frac{3}{2} h_{k-1} - \frac{1}{2} h_{k-d-1} < \frac{3}{2} h_{k-1},
$$
and because $h_1 =d/2$, $h_2 <3d/4$,
we arrive at
$$
h_k \leq \frac{ d}{2} \left(\frac{3}{2}\right)^{k-1}
$$
for $k \geq 1$. This inequality will be needed in the following.

So, let $\mu ( C ) = 2^l \ ( l \in \mathbb{N} )$.
If $C$ is already unimodular (i.e. $l=0$), we are done. If $C$ is not unimodular (i.e., $l\geq1$),
then  choose $i_1,\dots,i_m \in \{ 1, \ldots ,d \} $ with $ 1 \leq m \leq d$
and $i_j<i_k$ for $j<k$ such that the vectors $v_{i_1},\dots,v_{im}$ generate a minimal non-unimodular subcone. Then we set
$$
u = \frac{1}{2} ( v_{i_1} + \cdots + v_{i_m} ).
$$
That $u$ is an integral vector follows from a more general fact; see equation \eqref{vectorform} below.

Now we apply stellar subdivision to the cone $C$ by the vector $u$, which
gives us the cones
$$
C_{i_s} = \RR_+ v_1 + \cdots +\RR_+ v_{i_s-1}
+ \RR_+ u
+ \RR_+ v_{i_s+1} + \cdots + \RR_+ v_d,\quad 1 \leq s \leq m \leq d.
$$
\noindent For these cones of the first generation
(we regard the initial cone $C$ as the cone belonging to the $0$-th generation)
we have
$$
\mu ( C_{i_s}) =
\Big|\det \Bigl( v_1, \ldots, v_{i_{s}-1} , \frac{1}{2}
( v_{i_1} + \cdots +v_{i_m } ),v_{i_{s}+1}, \ldots ,v_d\Bigr)\Big| =
\frac{1}{2} \mu ( C) = 2^{l-1}.
$$
\noindent If $\mu (C_{i_s}) =1$, then the procedure stops. Otherwise it is
continued until we end with a triangulation of the initial cone $C$ by
unimodular cones of the $l$-th generation.

For the vectors $w_k$ which have been used for the stellar subdivisions
of the cones of the $(k-1)$-th generation
we get
$$
w_k \in h_k \Delta_C.
$$
\noindent We will prove this claim by induction on $k$.
For $k = 1$ it is obvious because
$( v_{i_1} + \cdots +v_{i_m})/2\in (d/2)\Delta_C $.
For $k > 1$, all generators $u_1,\ldots , u_d$ of a certain cone
$C'= \RR_+ u_1 + \cdots + \RR_+ u_d$ of the $(k-1)$-th generation
either belong to the initial vectors $v_1,\ldots , v_d$ or are vectors which
have been used for stellar subdivisions of cones of \emph{different} generations.
So by induction it follows
$$
u_i \in h_{n_i} \Delta_C, \quad n_i \leq k-1,
$$
\noindent where the $n_i$ are pairwise different.
The equality
$$
w_k =\frac{1}{2} ( u_{j_1} + \cdots +u_{j_v} ), \quad 1\leq v\leq d,
$$
\noindent immediately leads us to
$$
w_k \in \frac{1}{2} (h_{k-1} + \cdots +h_{k-d} ) \Delta_C = h_k \Delta_C,
$$
\noindent because the $h_i$ are increasing. Hence, we are done. 
\end{proof}

The theorem motivated us to come up with a triangulation algorithm which first triangulates
the underlying cone into subcones $D$ with $\mu(D) = 2^l$ $(l \in \mathbb{N})$. Such triangulations play a central role in this paper and therefore we give them a special name.

\begin{definition}
A \emph{$2$-triangulation} is a triangulation $\Sigma$ whose simplicial cones have multiplicities equal to powers of $2$.
\end{definition}

Since multiplicities of simplicial cones can be interpreted as orders of groups, our terminology is a close analogy to the notion of $2$-group.

The next, purely number theoretic lemma will be essential in the process of finding a $2$-triangulation of a given cone $C$. (By $\ld$ we denote the base $2$ logarithm.)

\begin{lemma}\label{odd}
Let $m$ and $p$ be two odd integers with $p/2< m < p$. Then there
exist natural numbers $s \leq \ld( p)$ and $t < p/2$ such that
$$
2 ^st = (2^{s-1} - 1 ) p + m.
$$
\end{lemma}

\begin{proof}
Because both $m$ and $p$ are odd, there exist a natural number $s>1$ and another odd number $q$ such that $p-m= 2^{s-1}q$. Now, let $t=(p-q)/2$. Then, $t$ is a natural number, since $p$ and $q$ are both odd. Hence,
$$
2^st=2^{s-1}(p-q)=(2^{s-1}-1)p + m,
$$
which proves the lemma.
\end{proof}

\begin{remark}
Improving Lemma \ref{odd} in the sense that we could find natural numbers
$s \ll \ld( p)$, $t < p/2$ and $x$ such that $2^s t = xp+m$ for given odd
numbers $m$ and $p$ with $p/2<m<p$ would critically affect the numerical quality of the bound we are going to give later on.
\end{remark}

In the following we will use an upper bound for the prime number counting function that
J. Rosser and L. Schoenfeld provided in \cite{RS62}.

\begin{theorem}\label{prime}
For $x > 0$ let $\pi (x)$ denote the number
of prime numbers $p$ with $p< x$. Then for all $x > 1$ we have
$$
\pi (x) < 1.25506 \cdot \frac {x}{\ln (x) }.
$$
\end{theorem}




As pointed out above, we want to subdivide a given simplicial cone $C$ into a $2$-triangulation by successive stellar subdivision. Therefore it is useful to replace a large prime number $p$ that divides $\mu(C)$ by smaller prime numbers in the passage from $C$ to the subcones resulting from a subdivision step.

Let the primitive vectors $v_1,\dots,v_d\in\ZZ^d$ generate a simplicial cone $C$ of dimension $d$, and let $U$ be the sublattice of $\ZZ^d$ spanned by these vectors. Then $\mu(C)$ is the index of $U$ in $\ZZ^d$, and each residue class has a representative in $$
\parf(v_1,\dots,v_d)=\{q_1v_1+\dots+q_dv_d: 0\le q_i < 1\}.
$$
If $p$ divides $\mu(C)$, then there is an element of order $p$ in $\ZZ^d/U$, and consequently there exists a vector
\begin{equation}\label{vectorform}
x= \frac{1}{p}\sum_{i=1}^d z_iv_i \in \ZZ^d, \qquad z_i\in \ZZ, \ 0\le z_i < p.
\end{equation}
The next lemma shows that we can find an element $x$ such that the coefficients $z_i$ avoid prime numbers between $3$ and $p$ for a subset of the indices $i$ whose size can be bounded by a function of $p$.

\begin{lemma}\label{avoid}
With the notation introduced, let $M\subset \{1,\dots,d\}$ such that 
$$
|M| \leq \frac{\ln(p)}{\tau},\quad \tau=1.25506.
$$
Then there exists an element $x$ of order $p$ modulo $U$ such that none of the coefficients $z_i$, $i\in M$, is an odd prime $<p$. 
\end{lemma}

\begin{proof} Let $b \rem a$ denote the remainder of $b$ modulo $a\neq 0$, chosen between $0$ and $|a|-1$.

Let $x$ be an element of order $p$ modulo $U$ given as above. Then all the elements
$$
x_j=\frac{1}{p}\sum_{i=1}^d (jz_i \rem p)v_i, \qquad j=1,\dots,p-1,
$$
lie in $\parf(v_1,\dots,v_d)$ and have order $p$ modulo $U$ since $p$ is prime.

We consider the maps $a\mapsto ja\rem p$, $0<j<p$. We must find a factor $j$ such that $jz_i\rem p$ is not an odd prime $<p$ for $i\in M$. If $z_i=0$, then this condition does not exclude any factor $j$. Otherwise it excludes $q$ factors where $q$ is the number of odd primes $<p$. In total we must exclude at most $|M|q$ factors $j$. But in view of Theorem \ref{prime} we have
$$
|M|q=|M|\left(\pi(p)-1\right)< |M|\left(\tau \cdot \frac{p}{\ln(p)}-1\right)\leq p-|M|\leq p-1
$$
for $|M|>0$. Furthermore, the lemma is obviously true for $|M|=0$.
\end{proof}


\section{The algorithm}

Before we describe the triangulation procedure precisely, we give an informal outline. The aim of this procedure is a $2$-triangulation $\Sigma$ of the original cone $C$ so that the generators of the cones $D\in\Sigma$ are relatively short with respect to the simplex $\Delta_C$.

For this purpose, we successively apply stellar
subdivision with carefully selected vectors $x \in C$ to the cone
$C = \RR_+ v_1+ \cdots + \RR_+ v_d \subset \RR^d$.
If $p$ is a prime divisor of $\mu(C)$, then there exists a vector
$$
x = \frac{1}{p}\left(\sum_{j=1}^d z_j v_j\right) \in \parf(v_1, \ldots, v_d) \setminus\{0\}.
$$
In order to end up with a $2$-triangulation of $C$, we want $z_j$ to be either a composite number or small, namely $z_j \leq p/2$. In general, $z_j$ cannot be expected to have this property. Therefore we add a certain multiple $kv_j$ $(k \in \mathbb{N})$ of $v_j$ to the vector $x$ if $z_j$ is a prime number and $z_j > p/2$. This results in a vector $x' \in C$
with
$$
x' = \frac{1}{p}\left(\sum_{j=1}^d z_j' v_j\right)
$$
such that all $z_j'$ are of the form $z_j'= 2^{g_j} t_j$ where $t_j\leq p/2$ or
$t_j$ is a composite number. By Lemma \ref{odd} we can achieve this goal.

Of course, we wish the vectors $x'$ to be as short as possible. We must
avoid the situation that both $z_j'$ is big and $v_j$ is a long vector because
then $x'$ would be long. Here the upper bound for the prime number counting
function comes into play via Lemma \ref{avoid} (see Lemma \ref{stellarvector} for a more detailed explanation), which guarantees that certain vectors are being multiplied by numbers $z_j'/p$ with $z_j'<p$.
\vspace{5 mm}


\noindent \hrulefill

\noindent\textbf{Power $2$ triangulation -- P2T}
\vspace{-3 mm}

\noindent \hrulefill

\begin{algorithmic}[1]

\REQUIRE The initial cone $C$

\ENSURE The 2-triangulation $\hat{T}(C)$ of $C$

\STATE $\hat{T}(C):= \{C\}$

\STATE $\hat{A}(C):=\{C\}$

\STATE $\xi_C(-i) := v_i$ for $i =1, \ldots , d$

\STATE $\tau :=1.25506$

\STATE $\xi_C(i) := 0$ for $i \in \mathbb{N}_0$

\WHILE {$\hat{T}(C)$ contains a cone $D= \RR_+ \xi_D(i_1)+ \cdots + \RR_+ \xi_D(i_d)$ (where $i_1 >i_2> \ldots >i_d\geq -d$) such that $\mu(D)$ is not a power of 2}

\STATE $p := \max\{p \in \mathbb{P}: \, p \mid \mu(D)\}$

\STATE {\textbf{FIND} $x = \nicefrac{1}{p}\left(\sum_{j=1}^d z_j \xi_D(i_j)\right)\in \parf(\xi_D(i_1), \ldots, \xi_D(i_d))\setminus\{0\}$ (which exists due to Lemma \ref{avoid} ) for which \\
\hspace{4 mm}\textbf{(1)} $z_j \notin \mathbb{P}$ or \\
\hspace{4 mm}\textbf{(2)} $z_j \leq p/2$ or \\
\hspace{4 mm}\textbf{(3)} $z_j = 2$ and $p=3$}\\
\quad for all $j \leq \ln(p)/\tau$

\FORALL{$j = \lfloor \ln(p)/\tau \rfloor +1 \ldots, d$}

\IF{$z_j \notin \mathbb{P}$ or $z_j \leq p/2$ or $z_j = 2$}

\STATE $z_j':= z_j$

\ELSE

\STATE $z_j':= z_j + k p$ with $k \in \mathbb{N}$ such that
$ z_j + k p = 2 ^st$ where $s \leq \ld( p)$ and
$t < p/2$ ( apply Lemma \ref{odd})

\ENDIF

\ENDFOR

\STATE $x':= \nicefrac{1}{p}\left(\sum_{j=1}^d z_j'w_j\right)$

\FORALL {$E \in \hat{T}(C)$ with $x' \in E$}

\STATE Apply stellar subdivision to $E$ by $x'$ (let
$E_j$ $(j= 1, \ldots, m)$ be the resulting cones)

\STATE $\hat{T}(C):= (\hat{T}(C) \setminus \{E\}) \cup \{E_j : \, j= 1, \ldots, m\}$

\STATE $\hat{A}(C):= \hat{A}(C) \cup \{E_j : \, j= 1, \ldots, m\}$

\ENDFOR

\STATE $\nu:= \max\{i : \xi_E (i)\neq 0\}$

\FORALL{$j = 1, \ldots, m$}

\FORALL{$k \leq \nu$}

\STATE $\xi_{E_j}(k) := \xi_{E}(k)$

\ENDFOR

\STATE $\xi_{E_j}(\nu +1) :=x'$

\ENDFOR

\ENDWHILE

\STATE Return $\hat{T}(C)$
\end{algorithmic}

\vspace{-3 mm}
\noindent \hrulefill
\vspace{5 mm}


The set $\hat{A}(C)$ contains the original cone $C$ and all cones being created in the course of the P2T algorithm. The set $\hat{T}(C)$ is a strict subset of $\hat{A}(C)$ unless $\mu(C)$ is a power of 2.   $\hat{A}(C)$ has been introduced out of technical reasons; it will help us to analyze certain properties of the resulting triangulation.

The line 8 in the P2T algorithm could be easily elaborated in a way how to determine the vector $x$ constructively, but because we are -- in the context of this article -- only interested in the existence of $x$, we leave it at this point.

The triangulation of the initial cone $C$ resulting from P2T
has the desired properties. We verify them in the following.

\begin{lemma}\label{stellarvector}
Let $D = \RR_+ w_1+ \cdots + \RR_+ w_d$ be a cone
to which we apply a stellar subdivision by a vector $x'$ in the P2T algorithm. Then
$x'$ is of the form
$$
x' = \frac{1}{p_{\max}}\left(\sum_{j=1}^d z_j' w_j\right), \quad
p_{\max} = \max\{p \in \mathbb{P}: \, p \mid \mu(D)\},
$$
such that:
\begin{enumerate}
\item for all $j$ we have $z_j' =2^{g_j} m_j$ with $g_j \in \mathbb{N}$, $g_j \leq \ld(p_{\max})$ and
\begin{enumerate}
\item $m_j \leq 2p_{\max}/3$ or
\item $m_j < p_{\max}$ is a composite number;
\end{enumerate}
\item 
\begin{enumerate}
\item $z_j'/p < 1$ for $j\le \ln(p)/\tau$,
\item $z_j'/p\le p/2$ for the remaining $j$.
\end{enumerate}
\end{enumerate}

\end{lemma}

This lemma has already been stated implicitly in lines 8 and 13 of the P2T algorithm.

The next definition will be
helpful in showing that the multiplicities of the cones in the final set
$\hat{T}(C)$ are relatively small and that the length of every chain of cones
$$
E_0=D \subset E_1 \subset E_2 \ldots \subset E_L=C,
$$ 
where $E_i$ is generated from $E_{i+1}$ by stellar subdivision and $D$ belongs to the resulting 2-triangulation of $C$, is relatively short.

\begin{definition}
Let $n$ be a natural number, $n= \prod_{i=1}^{\infty} p_i^{\alpha_i}$ be its prime
decomposition.
Then we define $\phi(n)= 2\left(\ld(n)- \eta(n)\right)$, where $\eta(n) = \sum_{i=1}^{\infty} \alpha_i$. (Hence $\phi(n) = \sum_{i=1}^{\infty}\alpha_i\left(2\ld(p_i)-2\right)$.)
\end{definition}

The function $\phi$ has some obvious nice properties, which we will need in the following.

\begin{lemma}\label{phi}\leavevmode
\begin{enumerate}
\item $\phi(ab) = \phi(a) + \phi(b)$ for $a,b \in \mathbb{N}$,
\item $\phi(a/b) = \phi(a) - \phi(b)$ for $a,b \in \mathbb{N}$, $b\mid a$,
\item $n= 2^s$ with $s\in\mathbb{N}$ if and only if $\lfloor\phi(n)\rfloor = 0$.
\end{enumerate}
\end{lemma}



\begin{lemma}\label{phimu}
Let $D,E \in \hat{A}(C)$ such that $E$ results from $D$ by stellar subdivision in the course of the P2T algorithm. Then
$$
\phi(\mu(E)) \leq \phi(\mu(D))-1.
$$
\end{lemma}

\begin{proof}
Due to lines 8 and 13 of the algorithm,
$$
\mu(E) = \mu(D)\frac{f}{p_{\max}} \cdot 2^l,
$$
where $p_{\max}=\max\{p\in \mathbb{P}:p|\mu(D)\}$ and $f,l \in \mathbb{N}$. Furthermore, $f$ is
\begin{enumerate}
\item either composite -- i.e., $f=u\cdot v <p_{\max}$ (with $u, v \in \mathbb{N}$) -- or
\item $f \leq 2 p_{\max}/3$.
\end{enumerate}
By Lemma \ref{phi} and because by definition $p_{\max}\mid\mu(D)$, we have
$$
\phi(\mu(E)) = \phi(\mu(D))-\phi(p_{\max}) + \phi(f)+ \phi(2^l) = \phi(\mu(D)) + \phi(f)-\phi(p_{\max}).
$$
In case (1)
$$
\phi(f)-\phi(p_{\max}) = \phi(u)+\phi(v) - 2\ld(p_{\max})-2 \leq -2,
$$ 
which proves the lemma in this case. In case (2) 
$$
\phi(f)-\phi(p_{\max}) \leq 2\cdot (\ld(p_{\max}) +\ld(2)-\ld(3)-\ld(p_{\max}))\leq -1,
$$ 
which proves the lemma for the second case.
\end{proof}

\begin{theorem}\label{ppttriang}
For a simplicial $d$-cone $C$ the P2T algorithm computes a $2$-triangulation of
$C$.
\end{theorem}

\begin{proof} The algorithm applies successive stellar subdivisions to the initial cone $C$. It stops when all multiplicities are powers of $2$, and that it stops after finitely many iterations follows from Lemma \ref{phimu}.
\end{proof}

\section{Bounds}

\begin{lemma}\label{chicone}
Let $D\in \hat{A}(C)$ be an arbitrary cone resulting from the P2T algorithm. Furthermore, we define
$$
\chi(D)= \max \{i: \xi_D(i) \neq 0\}.
$$
Then
$$
\chi(D) \leq \phi(\mu(C))-1.
$$
\end{lemma}

\begin{proof}
Let $D\in \hat{A}(C)$. By the algorithm, there is chain of cones
$$
E_0=D \subset E_1 \subset E_2 \ldots \subset E_L=C
$$ 
such that $E_i$ is generated from $E_{i+1}$ by stellar subdivision. Lemma \ref{phimu} implies that $\phi(\mu(D)) \leq \phi(\mu(C))-L$. On the other hand, by construction, $\chi(D) = \chi(C)+L$, where $\chi(C) =-1$. Therefore
$$
\chi(D) =L-1\leq \phi(\mu(C))-\phi(\mu(D))-1.
$$ 
This proves the lemma.
\end{proof}

\begin{theorem}\label{mubound}
For all $D \in \hat{A}(C)$ we have
$$
\mu (D) \leq 2^{1/2\cdot \ld(\mu(C))\cdot (\ld(\mu(C))+3)}
$$
\end{theorem}

\begin{proof}
By the algorithm, there is a chain of cones 
$$
D=E_0 \subset E_1 \subset E_2 \ldots \subset E_{L}=C,\qquad L \in \mathbb{N}_0
$$ 
such that $E_i$ is generated from $E_{i+1}$ by stellar subdivision.
Furthermore, let $p_{\max}(n)= \max\{p\in \mathbb{P}:p\mid n\}$ for any natural number $n$. 
Then obviously $p_{\max}(\mu(E_{i+1}))\geq p_{\max}(\mu(E_{i}))$ (see lines 8 and 13 of the algorithm).

Now, choose $s$ such that $p_{\max}(\mu(E_{i}))\leq 3$ for all $i\leq s$. Due to lines 8 and 13 of the algorithm, we have 
$$
\mu(E_{i}) = \frac{z_i}{p} \cdot\mu(E_{i+1})\leq \mu(E_{i+1}),\qquad 0\leq i< s
$$
since $z_i \leq p$ and $p\in \{2,3\}$.

On the other hand, again by lines 8 and 13, we have that
$$
\mu(E_{i}) = \mu(E_{i+1})\frac{f}{p_{\max}(\mu(E_{i+1}))}\cdot 2^l,\qquad i >s,
$$
where $f,l \in \mathbb{N}$. Furthermore, 
\begin{enumerate}
\item $f$ is either composite -- i.e., $f=u\cdot v <p_{\max}$ (with $u, v \in \mathbb{N}$) -- (see lines 8 and 13 of the algorithm) and $l=0$, or 
\item $f\leq p/2$ and $l \leq \ld(p_{\max}(\mu(E_{i+1})))$.
\end{enumerate}
In both cases it follows that 
$$
\frac{\mu(E_{i})}{2^{\eta(\mu(E_{i}))}} \leq
\frac12 \cdot \frac{\mu(E_{i+1})}{2^{\eta(\mu(E_{i+1}))}},\qquad i>s.
$$

On the other hand, $p_{\max}(n) \leq 2n/2^{\eta(n)}$ for every natural number $n$. 
Therefore, 
$$
p_{\max}(\mu(E_{i}))\leq \frac{2\mu(E_{i})}{2^{\eta(\mu(E_{i}))}} \leq \frac{1}{2^{L-i}} \cdot \mu(C),\qquad i>s
$$
This implies that $L-s\leq \lfloor\ld(\mu(C))\rfloor$. Otherwise we would have that $p_{\max}(\mu(E_{s+1}))\leq 2$.
Furthermore, it follows that
$$
\mu(E_i)\leq p_{\max}(\mu(E_{i+1})) \cdot \mu(E_{i+1})\leq 
\frac{1}{2^{L-i-1}} \cdot \mu(C) \cdot \mu(E_{i+1}),\qquad i>s.
$$
For $\ld(\mu(C))\geq 2$ we have that
\begin{multline*}
\mu(D)=\prod_{i=0}^s \frac{\mu(E_i)}{\mu(E_{i+1})} \cdot \prod_{i=s+1}^{L-1}\frac{\mu(E_i)}{\mu(E_{i+1})} \cdot \mu(E_L) \leq 
\mu(C) \prod_{i=s+1}^{\lfloor\ld(\mu(C))\rfloor+s-1} \frac{\mu(C)}{2^{L-i-1}}\\ \leq 
2^{1/2\cdot \ld(\mu(C))\cdot (\ld(\mu(C))+3)},
\end{multline*}
For $\mu(C)=3$ the algorithm stops after the first iteration, because there is a vector $x$ as given in line 8 of the algorithm, where $z_j\in \{0,1,2\}$ for all $j$. Hence, the resulting cones do have multiplicities equal to 1 or 2. Therefore, for all cones $D$ we have that 
$$
\mu(D)\leq 2\leq 2^{1/2\cdot \ld(\mu(C))\cdot (\ld(\mu(C))+3)}.
$$
Furthermore, for $\mu(C)\in \{1,2\}$ the algorithm even stops before the first iteration (see lines 1 and 6), which implies that $C=D$. Hence, in this case
$$
\mu(D)=\mu(C)\leq 2^{1/2\cdot \ld(\mu(C))\cdot (\ld(\mu(C))+3)}
$$
which finishes the proof.
\end{proof}


The next theorem is the central numerical consequence
resulting from the P2T algorithm, namely a length bound on the vectors involved. In the theorem and its proof we use the notation of P2T.

\begin{theorem}\label{vectorlength}
Let $D \in \hat{T}(C)$. Then, for all $s\geq 0$:
$$
\xi_D(s) \in \biggl( \frac{d}{2} \cdot \mu(C) \cdot 4^s \biggr)\Delta_C.
$$
\end{theorem}

\begin{proof}
To simplify notation we set $\mu=\mu(C)$. We prove the theorem via induction on $s$. So, let $s =0$.
If $\xi_D(0) = 0$, there is nothing to prove. 

So suppose that $\xi_D(0) \neq 0$. By the construction of $\xi_D(0)$ it follows that this vector was used for the stellar subdivision of the initial cone $C$. Hence, $\xi_D(0)$ is of the form
$$
\xi_D(0) = \frac{1}{p}\sum_{i = 1}^d z_i' v_i \in \ZZ^d \setminus \{0\}.
$$
where
$z_i'/p \leq p/2 \le \mu/2$ for all $i$ (Lemma \ref{stellarvector}). Therefore $x \in (d/2) \mu\Delta_C$, 
which finishes the case $s=0$. 

For the induction step assume the statement is true for $s$ replaced by $s-1\geq 0$. Again there is nothing to prove if $\xi_D(s) = 0$. Otherwise $\xi_D(s) \neq 0$ is a vector used for stellar subdivision. With the same notation as above, it follows
by construction of $\xi_D(s)$ that
$$
\xi_D(s) = \frac{1}{p}\left(\sum_{i = 1}^d z_i' \xi_D(j_i) \right)\in \ZZ^d \setminus \{0\}
$$
such that $s> j_1>j_2> \ldots > j_d$. Furthermore $p$ is a prime number $\le \mu$. Now, let $l$ be chosen such that $j_{l} > -1 \geq j_{l+1}$ (or $l=d$), implying that $\xi_D(j_l+1), \ldots, \xi_D(j_d)\in \{v_1, \ldots, v_d\}$.
We set
$$
q=\Bigl\lfloor \frac{\ln(p)}{\tau}\Bigr\rfloor
$$
and distinguish three cases.

(1) $q=0$. This is equivalent to $p=3$. It follows by induction that
$$
\xi_D(s) \in \left(\sum_{i=1}^{l} \left(\frac{d}{2}\mu 4^{j_i}\right)\frac{z_i'}{3} +
\sum_{i=l+1}^{d}\frac{z_i'}{3}\right) \Delta_C,
$$
because $s> j_1>j_2> \dots> j_d$. But, for $p=3$ we have that $z_i'<3$ for all $i$ (see line 8 of the algorithm). Hence,
$$
\xi_D(s) \in \frac{d}{2} \mu\biggl( \sum_{i=1}^{l} 4^{j_i}+1\biggr)\Delta_C\subset \left(\frac{d}{2} \mu 4^{s}\right)\Delta_C,
$$
since $\mu \geq p=3$.

(2) $l\leq q$ (and $q\neq0$). Then, again, it follows by induction that
$$
\xi_D(s) \in \left(\sum_{i=1}^{l} \left(\frac{d}{2}\mu 4^{j_i}\right)\frac{z_i'}{p} +
\sum_{i=l+1}^{d}\frac{p}{2}\right) \Delta_C,
$$
In the algorithm $z_i'/p<1$ for all $i\leq \ln(p)/\tau$ and therefore $z_i'/p<1$ for all $i\leq l$ . Hence,
$$
\xi_D(s) \in \frac{d}{2} \mu\biggl( \sum_{i=1}^{l} 4^{j_i}+1\biggr)\Delta_C\subset \left(\frac{d}{2} \mu 4^{s}\right)\Delta_C,
$$
which finishes the argument in this case.

(3) $l>q>0$. By induction it follows that
$$
\xi_D(s) \in \left(\sum_{i=1}^{q} \left(\frac{d}{2}\mu 4^{j_i}\right)\frac{z_i'}{p}
+ \sum_{i=q +1}^{l} \left( \frac{d}{2}\mu 4^{j_i}\right)\frac{z_i'}{p}
+ \sum_{i=l+1}^{d} \frac{p}{2} \right) \Delta_C.
$$
From the first two sums we can extract the factor $(d/2)\mu$ and bound the third summand by $(d/2)\mu$.

Because $s> j_1>j_2> \ldots > j_d$ and $z_i'/p<1$ for all $i\leq q$, as well as $z_i'/p<p/2$ for $i> q$, we have that
$$
\sum_{i=1}^{q} 4^{j_i}\frac{z_i'}{p} +
\sum_{i=q +1}^{l} 4^{j_i}\frac{z_i'}{p}+1
\leq
\sum_{i=1}^{q} 4^{j_i} +1 +
\sum_{i=q +1}^{l} 4^{j_i}\frac{p}{2}.
$$
Furthermore, $\sum_{i=k}^{l} \lambda^{i}= \frac{\lambda^{l+1}-\lambda^{k}}{\lambda-1}$ for each $\lambda \neq 1$ and $k,l \in \mathbb{N} \ (k\leq l)$. Hence,
$$
\sum_{i=1}^{q} 4^{j_i} +1\leq \frac{1}{3}\cdot (4^{j_1+1}-4^{j_q})+1\leq \frac{1}{3}\cdot 4^{j_1+1},
$$
because $l>q$ and $j_l\geq 0$ which implies that $j_q \geq 1$. Therefore,
$$
\sum_{i=1}^{q} 4^{j_i} +1+
\sum_{i=q +1}^{l} 4^{j_i}(p/2)\leq
\frac{1}{3} \cdot 4^{j_1+1} + \frac{p}{6} \cdot 4^{j_{q+1}+1}
$$
Note that $j_i + r \leq j_{i-\lfloor r \rfloor}+1$ for $r\in \RR_+$. 
It implies
$$
j_{q+1} = j_{q+1}+q-q\leq j_1-q,
$$ 
hence
$$
\frac{1}{3} \cdot 4^{j_1+1} + \frac{p}{6} \cdot 4^{j_{q+1}+1}\leq 
\frac{1}{3} \cdot 4^{j_1+1} + \frac{p}{6} \cdot 4^{j_1-q+1}
$$
On the other hand,
$$
4^q=4^{\lfloor \ln(p)/\tau\rfloor}\geq 4^{\ln(p)/\tau-1}\geq \frac{p^{\ln(4)/\tau}}{4}\geq \frac{p}{4}.
$$
Finally,
$$
\frac{1}{3} \cdot 4^{j_1+1} + \frac{p}{6} \cdot 4^{j_1-q+1}\leq \frac{1}{3} \cdot 4^{j_1+1} +\frac{4}{6}\cdot 4^{j_1+1}\leq 4^s,
$$
which finishes the proof.
\end{proof}
\begin{theorem} \label{vectorlength2} Let $D \in \hat{T}(C)$. Then, for all $s\in \mathbb{Z}$ and $d\geq 2$:
$$
\xi_D(s) \in \biggl( \frac{d}{2} \cdot \mu(C) \cdot 4^{\phi(\mu(C))}\biggr)\Delta_C.
$$
\end{theorem}
\begin{proof} The theorem follows from Lemma \ref{chicone} and theorem \ref{vectorlength}. Due to Theorem \ref{vectorlength} 
$$
\xi_D(s) \in \biggl( \frac{d}{2} \cdot \mu(C) \cdot 4^s \biggr)\Delta_C
$$ 
for $s\geq 0$, but on the other hand $\max \{i: \xi_D(i) \neq 0\} \leq \phi(\mu(C))-1$ by Lemma \ref{chicone}, which shows that the theorem is true for all $s\geq 0$.
Furthermore, $\xi_D(s) \in \Delta_C$ for $s< 0$ by definition. Because $d\geq 2$ and $\phi(n)\geq 0$ for all natural numbers $n\geq 1$ one has
$$
\xi_D(s) \in \Delta_C \subset \biggl( \frac{d}{2} \cdot \mu(C) \cdot 4^{\phi(\mu(C))} \biggr)\Delta_C
$$
for all $s < 0$, which finishes the proof.

\end{proof}

If we now collect the results from Theorem \ref{mubound}, Theorem \ref{vectorlength2} and Lemma \ref{power2} and additionally keep
in mind that the P2T algorithm produces a $2$-triangulation of the cone $C$, then we arrive
at the desired result.

\begin{theorem}\label{resultphi}
Every simplicial $d$-cone
$C = \RR_+ v_1 + \cdots + \RR_+ v_d \subset \RR^d , \, d\geq 2$, has
a unimodular triangulation $C = D_1 \cup \ldots \cup D_t$ such that for all $i$
$$
\Hilb(D_i) \subset \left( \frac{d^2}{4} \cdot \mu(C)\cdot 4^{\phi(\mu(C))} \cdot \left(\frac{3}{2}\right)^{1/2\cdot \ld(\mu(C))\cdot (\ld(\mu(C))+3)}\right)
\Delta_C.
$$
\end{theorem} 


Using an upper bound for the function $\phi$, we can simplify the bound somewhat:

\begin{corollary}\label{simpler}
Let $\varepsilon= 5+3/2\cdot \ld(3/2)$ and $\varrho = 1/2 \cdot \ld(3/2)$. So, $\varepsilon \approx 5.88$ and $\varrho \approx 0.29$.
Then every simplicial $d$-cone
$C = \RR_+ v_1 + \cdots + \RR_+ v_d \subset \RR^d , \, d\geq 2$, which is not already unimodular (i.e., $\mu(C)>1$) has a unimodular triangulation $C = D_1 \cup \ldots \cup D_t$ such that for all $i$
$$
\Hilb(D_i) \subset \left( \frac{d^2}{64} \cdot \mu(C)^{\varrho \cdot \ld(\mu(C)) + \varepsilon} \right)
\Delta_C.
$$
\end{corollary}

\begin{proof} The corollary follows from Theorem \ref{resultphi} and the fact that $\phi(n)\leq 2\ld(n)-2$ for all natural numbers $n>1$. This means that $4^{\phi(\mu(C))} \leq 1/16 \cdot \mu(C)^4$ and $$
\left(\frac{3}{2}\right)^{1/2\cdot \ld(\mu(C))\cdot (\ld(\mu(C))+3)}=
\mu(C)^{1/2 \cdot \ld(3/2) \cdot \ld(\mu(C)) + 3/2 \cdot \ld(3/2)}.
$$
\end{proof}



\begin{remark}\label{general}
In every iteration of the P2T algorithm it is guaranteed that the set $\hat{T}(C)$ constitutes a triangulation of the initial cone $C$. Even more so, after every iteration the set $\hat{T}(C)$ is a refined triangulation of the previous set $\hat{T}(C)$ via successive stellar subdivisions.
In particular, the algorithm would also end up with a triangulation of a cone $C$ in case we start it with a triangulation $\hat{T}(C)=\{D_1, \ldots, D_N\}$. Furthermore, the resulting cones would also coincide on the boundary, because in every iteration of the P2T algorithm stellar subdivision with a vector $x$ is applied to every cone, which contains $x$ (see lines 17 and 18 of the algorithm).
Now, every cone $C$ can be triangulated into simplicial cones $C'$ generated by extreme rays of $C$. So, if we start the algorithm with $\hat{T}(C)=\{C'_1, \ldots, C'_N\}$, it would end up with a $2$-triangulation of $C$ which is a refinement of the start triangulation. Hence, the algorithm essentially works for an arbitrary cone. We can replace the basic simplex $\Delta_C$ by the convex hull $\Gamma$ of the origin and the extreme integral generators of $C$, and $\mu(C)$ by the lattice normalized volume of $\Gamma$.
\end{remark}


\end{document}